\documentclass[12pt]{article}

\usepackage{amsthm,amsmath,amssymb}
\usepackage{graphics}

\newtheorem{la}{Lemma}[section]

\newtheorem{thm}[la]{Theorem}
\newtheorem{cor}[la]{Corollary}

\newtheorem{example}{Example}
\newtheorem{conj}{Conjecture}

\newcommand{\bex}{\begin{example}}
\newcommand{\eex}{\end{example}}

\def\qed{\hspace*{\fill}\vrule height6pt width4pt depth0pt\medskip}
\def\d{\displaystyle}

\date{}

\title{\bf From a Consequence of Bertrand's Postulate to Hamilton Cycles}

\author{Hong-Bin Chen\thanks{Department of Applied Mathematics, Feng Chia University, Taichung 40724, Taiwan
{\tt (Email: andanchen@gmail.com)} Research supported by MOST 105-2115-M-035-006-MY2 and partially done while the first author visited the Alfr\'{e}d R\'{e}nyi Institute of Mathematics in 2017.}
\and
Hung-Lin Fu\thanks{Department of Applied Mathematics, National Chiao Tung University, Hsinchu 30050, Taiwan
{\tt (Email: hlfu@math.nctu.edu.tw)} Research supported by MOST 106-2115-M-009-008}
\and
Jun-Yi Guo\thanks{Department of Mathematics, National Taiwan Normal University, Taipei 11677, Taiwan
{\tt (Email: davidguo@ntnu.edu.tw)} Research supported by MOST 106-2115-M-003-007}
}

\begin{document}

\maketitle


\begin{abstract}
A consequence of Bertrand's postulate, proved by L. Greenfield and S. Greenfield in 1998, assures that the set of integers $\{1,2,\cdots, 2n\}$ can be partitioned into pairs so that the sum of each pair is a prime number for any positive integer $n$. Cutting through it from the angle of Graph Theory, this paper provides new insights into the problem. We conjecture a stronger statement that the set of integers $\{1,2,\cdots, 2n\}$ can be rearranged into a cycle so that the sum of any two adjacent integers is a prime number. Our main result is that this conjecture is true for infinitely many cases.

\end{abstract}

\section{Introduction}

 Primes are the collection that has been extensively studied in Number Theory. This paper is motivated from a result, by Greenfield and Greenfield in 1998 \cite{greenfield}, concerning prime numbers.

 \begin{thm}\label{greenfield}\cite{greenfield}\\
 The set of integers $\{1, 2, 3, \cdots, 2n\}$, $n \geq 1$, can be partitioned into pairs $\{a_i,b_i\}$ such that $a_i+b_i$ is prime for all $i=1, 2, \cdots, n$.
 \end{thm}
 Taking $n=10$ for an example, $\{1, 2, \cdots, 20\}$ can be partitioned into pairs so that the sum of each pair is a prime number, as shown in Figure 1.
\begin{figure}
\includegraphics{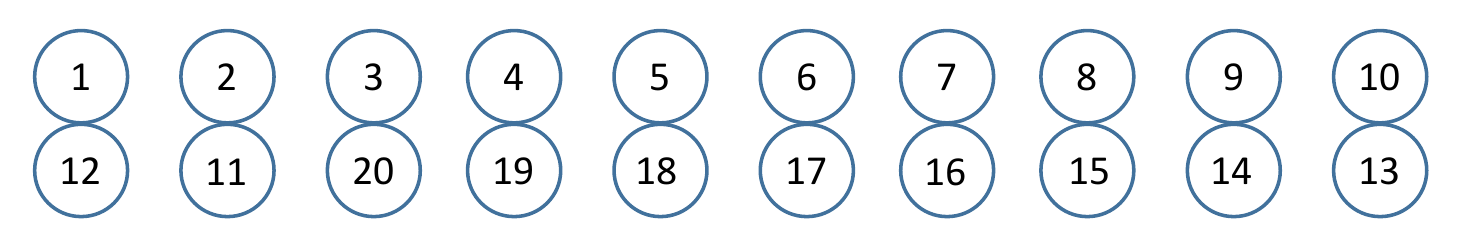}
\caption{$\{1, 2, \cdots, 20\}$ can be partitioned into prime pairs.}
\end{figure}

Theorem \ref{greenfield} was proved by L. Greenfield and S. Greenfield in 1998 \cite{greenfield} and reproduced by D. Galvin in 2006 \cite{Galvin}. This lovely result follows with an elegant proof from the well-known Bertrand’s Postulate, or sometimes called Bertrand-Chebyshev Theorem \cite{bertrand,tche}.

\begin{thm}\label{Bertrand}\cite{bertrand,tche}
For any positive integer $n>1$, there is at least a prime $p$ such that $n<p<2n$.
\end{thm}

As the proof of Theorem \ref{greenfield} is short, we demonstrate it in the following for self-contained. The proof is by induction on $n$. For $n=1$, it is trivial. Assume that the statement is true for any $k<n$. By Bertrand-Chebyshev Theorem, there exists a prime $p\in (2n,4n)$. We may assume that the prime $p=2n+k$ for some $1\leq k<2n$. Then the set of integers $\{k, k+1, \cdots, 2n-1, 2n\}$ obviously can be partitioned into pairs $\{k, 2n\}, \{k+1, 2n-1\}, \cdots$ and so on up to $\{n+\lfloor k/2\rfloor, n+\lceil k/2\rceil\}$ (the last is valid since $k$ is odd). Each of the pairs sums to the prime $2n+k$. By the induction hypothesis, the set of remaining integers $\{1, 2, \cdots, k-1\}$ has the property, too. Thus, the whole set of $\{1, 2, \cdots. 2n\}$ can be partitioned into pairs so that the sum of each pair is prime. This completes the proof.

From the Graph Theory point of view, it is natural to think of a graph that treats numbers as vertices and two vertices are adjacent if the sum of the corresponding numbers is a prime. Let's denote as $G_n$ such a graph of vertices $\{1, 2, \cdots, n\}$. Theorem \ref{greenfield} can be rephrased in the terminology of Graph Theory: a graph defined in this way has a perfect matching.
Inspired by Theorem \ref{greenfield}, we are interested in the structure of such a graph.

We first give a formal definition of the mentioned graph. For any positive integer $n$, define a graph $G_n=(V,E)$ with the vertex set $V=\{1,2,\cdots, n\}$ and $E=\{ij: i+j \text{ is prime}\}$. We call $G_n$ the {\it prime sum graph} of order $n$.

Expectedly, properties in $G_n$ can reflect some properties on prime numbers. Therefore, if one can provide some new insight from the angle of Graph Theory, then it could eventually raise awareness of new challenges of prime numbers.
Obviously, prime sum graphs are bipartite graphs, and thus are 2-colorable. This reflects a trivial fact that a prime (an edge) can be formed only by connecting an odd number to an even number. Each vertex in $G_n$ has degree $\frac{n}{\log n}$  as $n$ is sufficiently large. Thus, there are nearly $\frac{n^2}{2\log ⁡n}$ edges. This reflects a well-known result that the number of primes below $n$ is $\pi(n)\sim \frac{n}{\log ⁡ n}$.

Theorem \ref{greenfield} assures that $G_{2n}$ has a perfect matching. This result strongly relies on the existence of one prime $p$ between $2n$ and $4n$, as an immediate consequence of Bertrand-Chebyshev Theorem. In fact, there are some further extensions of Bertrand-Chebyshev Theorem that suggest the existence of more and more primes between $2n$ and $4n$ as $n$ is tending to infinity. Results alone this direction guarantee that every number connects to more than one numbers in a prime sum graph. Naturally, an interesting question is that \begin{center} ``can we say something more than a perfect matching hidden in the graph $G_{2n}$?'' \end{center} If there are two or more perfect matchings in $G_{2n}$, then it is likely that there exists a {\it Hamilton cycle}, a cycle that visits each vertex exactly once.

\begin{figure}
\begin{center}
\scalebox{0.6}{\includegraphics{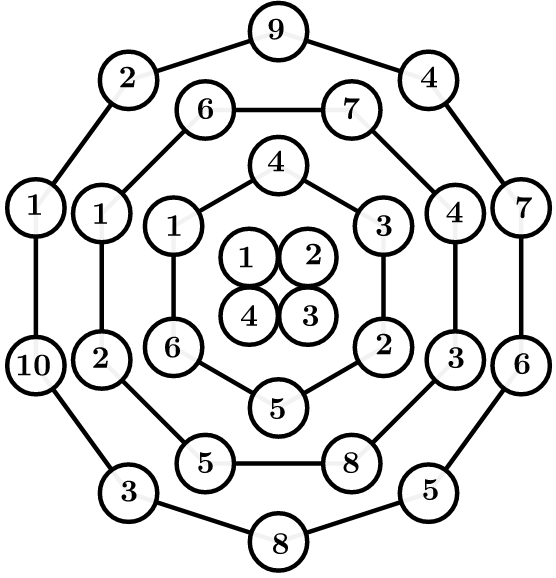}}
\caption{Examples for small $G_{2n}$'s with a Hamilton cycle.}
\end{center}
\end{figure}

In this paper, we pose an interesting conjecture concerning positive integers and prime numbers.
\begin{conj}\label{conjecture}
 The set of integers $\{1, 2, 3, \cdots, 2n\}$, $n \geq 2$, can be rearranged in a circle such that the sum of any two adjacent numbers is a prime. In other words, $G_{2n}$ contains a Hamilton cycle.
\end{conj}

We first prove the following result, which provides a sufficient condition for $G_{2n}$ being Hamiltonian.
\begin{thm}\label{main}
$G_{2n}$, $n\geq 2$, contains a Hamilton cycle if there exist two primes $p_1<p_2$ in $[1,2n]$ ($p_1$ can be 1) such that $2n+p_1$ and $2n+p_2$ are primes and gcd$\d\left(\frac{p_2-p_1}{2},n\right)=1$.
\end{thm}

For any positive integer $m$, let $E_m$ denote the set of the edges in a prime sum graph $G=(V,E)$ with the sum of two endpoints is $m$, i.e., $E_m=\{uv\in E(G): u+v=m\}$. Denote $\mathcal{P}$ as the set of all primes less than $4n$. A set $H\subseteq \mathcal{P}$ is called a {\it Hamiltonian prime set} of $G_{2n}$ if $\bigcup_{p\in H}E_p$ contains a Hamilton cycle in $G_{2n}$. As a consequence of Theorem \ref{main} by setting $p_1=1$ and $p_2=3$, we have the following corollary.

\begin{cor}\label{cor}
The set $\{3, 2n+1, 2n+3\}\subseteq \mathcal{P}$ is a Hamiltonian prime set of $G_{2n}$ for $n\geq 2$. In other words,
if $2n+1$ and $2n+3$ are twin primes, then $G_{2n}$ has a Hamilton cycle.
\end{cor}

It is worth mentioning that $\{3, 2n+1, 2n+3\}$ is the smallest Hamiltonian prime set of $G_{2n}$. The reason is that every vertex needs at least two edges to form a Hamilton cycle and thus the vertex $i$ requires at least two primes in its potential candidate set $[1+i, 2n+i]$ for each $1\leq i\leq 2n$.

On the one hand, Corollary \ref{cor} implies that if the well-known twin prime conjecture is true, then there are infinitely many $G_{2n}$'s that have a Hamilton cycle. On the other hand, the above discussion can also be interpreted as a big challenge of proving the existence of infinitely many Hamiltonian $G_{2n}$'s. The reason is that, to achieve this conclusion, it needs to prove that there are infinitely many prime triples (or quadruples) satisfying certain conditions. Many similar statements concerning prime numbers have been addressed and proposed in the literature. However, most are still unsolved and the twin prime conjecture is one of the mysteries.

A recent breakthrough in the twin prime conjecture is that Yitang Zhang showed in 2013 \cite{zhang} for the first time that we will never stop finding pairs of primes that are within a bounded distance — within 70 million. Soon after, dozens of outstanding researchers in the world work together to improve Zhang's 70 million bound, bringing it down to 246 \cite{maynard,polymath}.

\begin{thm}\label{246}\cite{maynard,polymath}
There are infinitely many pairs of primes $(p_i,p'_i)$ such that $p'_i-p_i\leq 246$.
\end{thm}

Thanks to the breakthrough, we can prove the following main result by combining Theorem \ref{main} and Theorem \ref{246} with an elaborate argument.

\begin{thm}\label{mainthm}
There are infinitely many $G_{2n}$'s that have a Hamilton cycle.
\end{thm}

This paper partly improves a Number Theory result by Greenfield and Greenfield \cite{greenfield} with new insights from Graph Theory. We pose Conjecture \ref{conjecture} and show that it is true for infinitely many cases. Although the main result is still far from our conjecture, the value of this paper is in calling awareness to the discovery of a simple, previously unnoticed property of numbers and its connection to several of the most central concepts in graphs and prime numbers.

The rest of this paper is organized as follows. Section 2 and Section 3 are the proofs of Theorems \ref{main} and \ref{mainthm}, respectively. Finally, we conclude this paper in Section 4.


\section{Proof of Theorem \ref{main}}

Consider a balanced bipartite graph $B=(X\cup Y,E)$ with vertex set $V=X\cup Y$, where $|X|=|Y|$, and edge set $E$ which
contains some edges with one vertex in $X$ and the other one in $Y$. Let $X=\{x_1,x_2,\cdots, x_n\}$ and $Y=\{y_1, y_2, \cdots, y_n\}$. For $0\leq i\leq n-1$, define {\it $i$-difference set} as the set consisting of edges whose indexes in $X$ and $Y$ differ by $i$ (mod $n$) and denote by $D_i=\{x_jy_k: i\equiv j-k ~(\text{mod}~ n)\}$. A result well known in Graph Theory assures as the following.

\begin{thm}\label{difference}
Let $B$ be a balanced bipartite graph of order $2n$, $n\geq 2$, as defined above and $s, t$ be two integers in $[0,n-1]$. Then $D_s\cup D_t$ forms a Hamilton cycle in $G$ if $\gcd(|t-s|,n)=1$.
\end{thm}

Recall that $E_m=\{uv\in E(G): u+v=m\}$. To prove the theorem, it suffices to show that $E_{p_1}\cup E_{2n+p_1}$ and $E_{p_2}\cup E_{2n+p_2}$ form two difference sets $D_s$ and $D_t$ for some $s$ and $t$, respectively, with $\gcd(|t-s|,n)=1$. Then, by Theorem \ref{difference}, $G_{2n}$ has a Hamilton cycle.

Consider $\{1,2,\cdots,2n\}=X\cup Y$ with $X=\{x_1,x_2,\cdots, x_n\}$ and $Y=\{y_1, y_2, \cdots, y_n\}$, where $x_j=2j-1$ for $1\leq j\leq n$ and $y_k=2n-2(k-1)$ for $1\leq k\leq n$. Then $D_{\frac{p_1-1}{2}}=\{x_{j}y_{k}: j-k\equiv \frac{p_1-1}{2} ~(\text{mod}~ n)\}$ and $D_{\frac{p_2-1}{2}}=\{x_{j}y_{k}: j-k\equiv \frac{p_2-1}{2} ~(\text{mod}~ n)\}$. Next, we first prove that $D_{\frac{p_1-1}{2}}=E_{p_1}\cup E_{2n+p_1}$. Assume that $x_jy_k\in D_{\frac{p_1-1}{2}}$.\\
 If $j\geq k$, then we have \begin{align*}
& j-k=\frac{p_1-1}{2}\\
      \Leftrightarrow ~& 2n+2(j-k)+1=2n+p_1\\
      \Leftrightarrow ~& \underbrace{(2j-1)}_{x_j}+\underbrace{2n-2(k-1)}_{y_k}=2n+p_1.
\end{align*} The last equality implies that $x_jy_k\in E_{2n+p_1}$, and vice versa. Thus, if $x_jy_k\in E_{2n+p_1}$ and $j\geq k$, then  $x_jy_k\in D_{\frac{p_1-1}{2}}$.\\ If $j<k$, then we have  \begin{align*}
& j-k=\frac{p_1-1}{2}-n\\
      \Leftrightarrow ~& 2n+2(j-k)+1=p_1\\
      \Leftrightarrow ~& \underbrace{(2j-1)}_{x_j}+\underbrace{2n-2(k-1)}_{y_k}=p_1.
\end{align*} This implies that $x_jy_k\in E_{p_1}$, and vice versa. Thus, if $x_jy_k\in E_{p_1}$ and $j< k$, then  $x_jy_k\in D_{\frac{p_1-1}{2}}$. Therefore, $D_{\frac{p_1-1}{2}}=E_{p_1}\cup E_{2n+p_1}$. A similar argument shows that $D_{\frac{p_2-1}{2}}=E_{p_2}\cup E_{2n+p_2}$.

It is easy to see that $\d \gcd\left(\frac{p_2-1}{2}-\frac{p_1-1}{2},n\right)=1$ as $\d \gcd\left(\frac{p_2-p_1}{2},n\right)=1$. By Theorem \ref{difference}, $D_{\frac{p_1-1}{2}}\cup D_{\frac{p_2-1}{2}}$ forms a Hamilton cycle in $G_{2n}$. \qed

\section{Proof of Theorem \ref{mainthm}}

Suppose that there exists a number $g$ such that there are infinitely many prime pairs $(p,p')$ satisfying the condition $p'-p=g$. Take $g=12$ for example.
\begin{itemize}
\item[(1)] Then, there must be infinitely many prime pairs $(p,p')$ with $p'-p=12$ and satisfying one of the following forms (all potential cases of a prime (mod 12)):
\[
\begin{tabular}{rcl}
$p$ & = & $12k+1$, \\
$p$ & = & $12k+5$, \\
$p$ & = & $12k+7$, \\
$p$ & = & $12k+11$. \\
\end{tabular}
\]
Notice that $12k+3$ and $12k+9$ are excluded because they are not primes.
\item[(2)] For each form, find an explicit representative as shown in Table 1.
  \begin{table}[h!]\label{tab1}
\[
\begin{tabular}{c|c}
Form & Representative \\
\hline
$p=12k+1$ & (1,13)\\
$p=12k+5$ & (5,17)\\
$p=12k+7$ & (7,19)\\
$p=12k+11$ & (11,23)\\
\end{tabular}
\]
\caption{Representatives}
  \end{table}

\item[(3)] Assume that it is exactly the form $p=12k'+1$ such that there are infinitely many prime pairs $(p,p')$ with $p'-p=12$. Then, according to Theorem \ref{main}, one can conclude that there are infinitely many $G_{2n}$'s having a Hamilton cycle by setting $p_1=11, p_2=23$, $2n+p_1=p$ and $2n+p_2=p'$.
    It is easy to verify that $2n=(12k'+1)-11=12k+2$ and thus $n=6k+1$ for some $k$. Therefore, we have $\d\gcd\left(\frac{p'-p}{2},n\right)=\gcd\left(6,6k+1\right)=1$, as desired.

    Similarly, for each of the other forms we can conclude the same by setting $p_1$ and $p_2$ properly, as in Table 2.
     \begin{table}[h!]\label{tab2}
      \[
\begin{tabular}{l|l|l}
Form & $(p_1,p_2)$ & gcd condition \\
\hline
$p=12k'+1$ & (11,23) & gcd$(6,6k+1)=1$ \\
$p=12k'+5$ & (7,19) &  gcd$(6,6k+5)=1$ \\
$p=12k'+7$ & (5,17) &  gcd$(6,6k+1)=1$ \\
$p=12k'+11$ & (1,13) &  gcd$(6,6k+5)=1$
\end{tabular}
\]
\caption{Match and gcd}
\end{table}

\item[(4)] It is known by Theorem \ref{246} \cite{maynard,polymath} that such a number $g\leq 246$ exists, but not knowing which number. We can check (1), (2), and (3) for all possibilities of $g=2,4,6,\cdots,246$ to see whether it is true or not. Since for each $g$ there are only finite cases, this process can be done simply by computers. As expected, for each $g$ and for each potential form $p=gk+t$, we successfully find a representative $(p_1,p_2)$ as a match with $(p,p')$ such that $\d \gcd(\frac{g}{2},\frac{g}{2}\cdot k+s)=1$.
    Because there are too many cases (total 6170 cases) to demonstrate in this paper, we list only the cases of $g=246$ below. The other cases can be found in the file of full lists in the webpage (http://www.davidguo.idv.tw/temp/Data.pdf).
\end{itemize}

The above discussion implies that there are infinitely many $G_{2n}$'s with a Hamilton cycle, verifying Conjecture 1 for infinitely many cases. The proof is complete.\qed

 \begin{table}[h!]\label{tab3}
\footnotesize
      \[
\begin{tabular}{c|l|l||c|l|l}
Form & $(p_1,p_2)$ & gcd condition & Form & $(p_1,p_2)$ & gcd condition\\
\hline
$246k'+1$ & (5,251) & gcd$(123,123k+121)=1$ &
$246k'+5$ & (2707,2953) & gcd$(123,123k+2)=1$\\
$246k'+7$ & (5,251) & gcd$(123,123k+1)=1$ &
$246k'+11$ & (2707,2953) & gcd$(123,123k+5)=1$\\
$246k'+13$ & (5,251) & gcd$(123,123k+4)=1$ &
$246k'+17$ & (2707,2953) & gcd$(123,123k+8)=1$\\
$246k'+19$ & (5,251) & gcd$(123,123k+7)=1$ &
$246k'+23$ & (2707,2953) & gcd$(123,123k+11)=1$\\
$246k'+25$ & (5,251) & gcd$(123,123k+10)=1$ &
$246k'+29$ & (2707,2953) & gcd$(123,123k+14)=1$\\
$246k'+31$ & (5,251) & gcd$(123,123k+13)=1$ &
$246k'+35$ & (2707,2953) & gcd$(123,123k+17)=1$\\
$246k'+37$ & (5,251) & gcd$(123,123k+16)=1$&
$246k'+43$ & (5,251) & gcd$(123,123k+19)=1$\\
$246k'+47$ & (2707,2953) & gcd$(123,123k+23)=1$&
$246k'+49$ & (5,251) & gcd$(123,123k+22)=1$\\
$246k'+53$ & (2707,2953) & gcd$(123,123k+26)=1$&
$246k'+55$ & (5,251) & gcd$(123,123k+25)=1$\\
$246k'+59$ & (2707,2953) & gcd$(123,123k+29)=1$&
$246k'+61$ & (5,251) & gcd$(123,123k+28)=1$\\
$246k'+65$ & (2707,2953) & gcd$(123,123k+32)=1$&
$246k'+67$ & (5,251) & gcd$(123,123k+31)=1$\\
$246k'+71$ & (2707,2953) & gcd$(123,123k+35)=1$&
$246k'+73$ & (5,251) & gcd$(123,123k+34)=1$\\
$246k'+77$ & (2707,2953) & gcd$(123,123k+38)=1$&
$246k'+79$ & (5,251) & gcd$(123,123k+37)=1$\\
$246k'+83$ & (991,1237) & gcd$(123,123k+38)=1$&
$246k'+85$ & (5,251) & gcd$(123,123k+40)=1$\\
$246k'+89$ & (2707,2953) & gcd$(123,123k+44)=1$&
$246k'+91$ & (5,251) & gcd$(123,123k+43)=1$\\
$246k'+95$ & (2707,2953) & gcd$(123,123k+47)=1$&
$246k'+97$ & (5,251) & gcd$(123,123k+46)=1$\\
$246k'+101$ & (2707,2953) & gcd$(123,123k+50)=1$&
$246k'+103$ & (5,251) & gcd$(123,123k+49)=1$\\
$246k'+107$ & (2707,2953) & gcd$(123,123k+53)=1$&
$246k'+109$ & (5,251) & gcd$(123,123k+52)=1$\\
$246k'+113$ & (2707,2953) & gcd$(123,123k+56)=1$&
$246k'+115$ & (5,251) & gcd$(123,123k+55)=1$\\
$246k'+119$ & (2707,2953) & gcd$(123,123k+59)=1$&
$246k'+121$ & (5,251) & gcd$(123,123k+58)=1$\\
$246k'+125$ & (2707,2953) & gcd$(123,123k+62)=1$&
$246k'+127$ & (5,251) & gcd$(123,123k+61)=1$\\
$246k'+131$ & (2707,2953) & gcd$(123,123k+65)=1$&
$246k'+133$ & (5,251) & gcd$(123,123k+64)=1$\\
$246k'+137$ & (2707,2953) & gcd$(123,123k+68)=1$&
$246k'+139$ & (5,251) & gcd$(123,123k+67)=1$\\
$246k'+143$ & (2707,2953) & gcd$(123,123k+71)=1$&
$246k'+145$ & (5,251) & gcd$(123,123k+70)=1$\\
$246k'+149$ & (2707,2953) & gcd$(123,123k+74)=1$&
$246k'+151$ & (5,251) & gcd$(123,123k+73)=1$\\
$246k'+155$ & (2707,2953) & gcd$(123,123k+77)=1$&
$246k'+157$ & (5,251) & gcd$(123,123k+76)=1$\\
$246k'+161$ & (2707,2953) & gcd$(123,123k+80)=1$&
$246k'+163$ & (5,251) & gcd$(123,123k+79)=1$\\
$246k'+167$ & (2707,2953) & gcd$(123,123k+83)=1$&
$246k'+169$ & (11,257) & gcd$(123,123k+79)=1$\\
$246k'+173$ & (2707,2953) & gcd$(123,123k+86)=1$&
$246k'+175$ & (5,251) & gcd$(123,123k+85)=1$\\
$246k'+179$ & (2707,2953) & gcd$(123,123k+89)=1$&
$246k'+181$ & (5,251) & gcd$(123,123k+88)=1$\\
$246k'+185$ & (2707,2953) & gcd$(123,123k+92)=1$&
$246k'+187$ & (5,251) & gcd$(123,123k+91)=1$\\
$246k'+191$ & (2707,2953) & gcd$(123,123k+95)=1$&
$246k'+193$ & (5,251) & gcd$(123,123k+94)=1$\\
$246k'+197$ & (2707,2953) & gcd$(123,123k+98)=1$&
$246k'+199$ & (5,251) & gcd$(123,123k+97)=1$\\
$246k'+203$ & (2707,2953) & gcd$(123,123k+101)=1$&
$246k'+209$ & (2707,2953) & gcd$(123,123k+104)=1$\\
$246k'+211$ & (5,251) & gcd$(123,123k+103)=1$&
$246k'+215$ & (2707,2953) & gcd$(123,123k+107)=1$\\
$246k'+217$ & (5,251) & gcd$(123,123k+106)=1$&
$246k'+221$ & (2707,2953) & gcd$(123,123k+110)=1$\\
$246k'+223$ & (5,251) & gcd$(123,123k+109)=1$&
$246k'+227$ & (2707,2953) & gcd$(123,123k+113)=1$\\
$246k'+229$ & (5,251) & gcd$(123,123k+112)=1$&
$246k'+233$ & (2707,2953) & gcd$(123,123k+116)=1$\\
$246k'+235$ & (5,251) & gcd$(123,123k+115)=1$&
$246k'+239$ & (2707,2953) & gcd$(123,123k+119)=1$\\
$246k'+241$ & (5,251) & gcd$(123,123k+118)=1$&
$246k'+245$ & (2707,2953) & gcd$(123,123k+122)=1$

\end{tabular}
\]
\caption{The cases of $g=246$.}
\end{table}

\section{Conclusions}

The notion of Hamilton cycles is one of the most central in modern Graph Theory, and many efforts have been devoted to obtain sufficient conditions for Hamiltonicity.
One of the oldest results is the theorem of Dirac \cite{dirac}, who showed that if the minimum degree of a graph $G$ on $n$ vertices is at least $\frac{n}{2}$, then $G$ contains a Hamilton cycle. This result is only one example of a vast
majority of known sufficient conditions for Hamiltonicity that mainly deal with fairly dense graphs. On the other hand, it appears that little has been known about Hamilton cycles in relatively sparse graphs. As mentioned previously, the minimum degree of $G_{2n}$ is roughly $\frac{n}{\log n} < \frac{n}{2}$. So, Dirac’s result or Dirac-type results for bipartite graphs \cite{erdos,moon} do not work in our case.

Things are different in random graphs. P\'{o}sa's result \cite{posa} shows that the probability that a random graph with $n$ vertices and $cn\log n$ edges, for sufficiently large $c$, contains a Hamilton cycle tends to 1 as $n$ tends to infinity. The number of edges in $G_{2n}$ is around $n^2/\log n$, which is more than $cn\log n$ in the P\'{o}sa's result. Although $G_{2n}$ is not a real random graph, connecting of the edges in $G_{2n}$ relies on prime numbers, which have been long thought to distribute `randomly' in some sense.

 It is worth mentioning that the sufficient condition in Theorem \ref{main} involves only four primes, and the number four is much less than $\pi(4n-1)$, the number of all candidates of primes in $G_{2n}$. Though the sufficient condition in Theorem \ref{main} is much stronger than Conjecture 1 in this sense, we believe that it is likely to be satisfied for any even number $2n\geq 4$. To support this point, we have verified the condition for all $2n<10^8$ by a computer search. A wealth of numerical evidence supports Conjecture \ref{conjecture}, but so far a proof has eluded us.


The numerical results of the existence of $p_2$ and  $2n+p_2$ in the sufficient condition also lead us to pose the following question.
\begin{center}
``Is there a prime $p<2n$ such that $2n+p$ is also a prime for any $2n\geq 4$?''
\end{center}
Obviously, the prime $2n+p$ is between $2n$ and $4n$. So, this can be viewed as a generalization of Bertrand's postulate on even numbers with an extra condition. In addition, the conjecture can also be interpreted as  a variant of Goldbach's conjecture that ``every even number $2n\geq 4$ is the difference of two primes $p$ and $2n+p$''. Chen's work \cite{chen} on Goldbach's conjecture also showed a related result that every even number is the difference between a prime and a product of two primes. However, it is still unknown ``whether every even number is the difference of two primes?''.  According to the discussion, it seems difficult to prove Conjecture \ref{conjecture} by showing directly that the sufficient condition holds for any even number $2n\geq 4$.


Another research direction is to consider an extension of the discussed problem alone the concept of {\it quasi-random graphs} by Chung, Graham and  Wilson \cite{chung}. Quasi-random graphs (or pseudo-random graphs) can be informally described as graphs whose edge distribution resembles closely that of a truly random graph $G(n,p)$ of the same edge density. The connection between pseudo-randomness and Hamiltonicity has been explored in several articles \cite{frieze2000,frieze2002,krivelevich}. Although the distribution of prime numbers in integers is not random, the prime sum condition on the discussed problem can be relaxed to any subset of integers, nothing to do with primes. By choosing these integers randomly, one can construct a new class of graphs that might behave like pseudo-random graphs. We believe that proving Hamiltonicity in the new class of graphs is an interesting and challenging task, and will lead to  more interesting problems.


\section*{Acknowledgments}
The first author would like to thank Professor Mikl\'{o}s Simonovits for valuable suggestions.

\end{document}